\documentclass{amsart}
\usepackage{enumerate}
\usepackage{amsthm,amscd,amssymb,verbatim,epsf,amsmath,amsfonts,mathrsfs,graphicx,enumitem,}

\usepackage[linktocpage,colorlinks=true,linkcolor=blue,citecolor=blue]{hyperref}
\usepackage{lineno}
\newcommand{\pard}[2]{\frac{\partial #1}{\partial #2}}
\newcommand{\ip}[2]{\left \langle #1 , #2 \right\rangle}
\newcommand{\n}{\nabla}
\newcommand{\ov}[1]{\overline{#1}}
\renewcommand{\tilde}[1]{\widetilde{#1}}
\newcommand{\bb}[1]{\mathbb{#1}}
\newcommand{\ra}{\rightarrow}
\renewcommand{\l}{\lambda}
\newtheorem{thm}{Theorem}
\newtheorem{lem}{Lemma}
\newtheorem{prop}{Proposition}

\newtheorem{rem}{Remark}
\newtheorem*{assump*}{Assumptions on $\psi$}

\title{ Prescribed $k$-symmetric curvature hypersurfaces in de Sitter space}

\author[D. Ballesteros-Ch\'avez]{Daniel Ballesteros-Ch\'avez}
\address{Daniel Ballesteros Ch\'avez\\
  Department of Mathematical Sciences\\
  University of Durham\\
  Stockton Road\\
  Durham\\
  DH1 3LE\\
  United Kingdom.}
\email{daniel.ballesteros-chavez@durham.ac.uk}

\author[W. Klingenberg]{Wilhelm Klingenberg}
\address{Wilhelm Klingenberg\\
 Department of Mathematical Sciences\\
 University of Durham\\
 Durham DH1 3LE\\
 United Kingdom.}
\email{wilhelm.klingenberg@durham.ac.uk }

\author[B. Lambert]{Ben Lambert}
\address{Ben Lambert\\Department of Mathematics, University College London, Gower Street, London, WC1E 6BT, United Kingdom}
\email{b.lambert@ucl.ac.uk}

\thanks{The first author was supported by CONACYT-Doctoral scholarship no.~411485, the third author was supported by a Leverhulme Trust Research Project Grant RPG-2016-174.}
\date{}
\subjclass[2010]{35J60, 53C50, 53C42}

\begin{document}
\maketitle

\begin{abstract}
We prove existence of compact spacelike hypersurfaces with prescribed k - curvature in de Sitter space, where the prescription function depends on both space and the tilt function.
\end{abstract}

\section{Introduction}
\noindent
We consider the existence problem for embedded compact spacelike hypersurfaces $\Sigma$ in de Sitter space $S^{n+1}_1$  satisfying a prescribed curvature equation of the form
\begin{equation}
  \label{ActualEquation}
   H_k^\frac 1 k(\lambda[A]) = \psi.
\end{equation}
\noindent
Here $1\leq k\leq n$ is fixed, $A$ is the second fundamental form of $\Sigma $, $\lambda[A]=(\lambda_1, \ldots, \lambda_n)$ are the eigenvalues of the shape operator $A^i_j$, and $H_k$ is the $k$-th normalised symmetric polynomial in $\lambda$, that is
\begin{equation}
    S_{k}(\lambda) :=  \sum_{1\leq i_1<\cdots<i_k\leq n}\lambda_{i_1}\cdots\lambda_{i_k}, \qquad\qquad H_k:= \begin{pmatrix} n \\ k \end{pmatrix}^{-1}S_k.
\end{equation}
\noindent
Furthermore, we will choose our prescription function  $\psi$ to depend on both position in $S^{n+1}_1$ and the tilt, $\tau$, which is a measure of how spacelike $\Sigma$ is at the position, see \eqref{Tiltdefinition} below. We note that, as every compact embedded spacelike surface in $S^{n+1}_1$ may be written as a graph, (\ref{ActualEquation}) may be rewritten as a fully nonlinear elliptic second order partial differential equation in the graph function. Our main result is that, assuming some natural structural assumptions on $\psi$, the prescribed curvature equation (\ref{ActualEquation}) has a smooth spacelike solution $\Sigma$. 

\noindent
The existence of solutions of such equations was studied in \cite{CNS3} by L. Caffarelli, L. Nirenberg and J. Spruck. In \cite{CNS4}, they proved the existence of starshaped hypersurfaces in Euclidean space
with prescribed $k$-symmetric curvature using an priori $C^{2,\alpha}$ estimate needed to carry out
the continuity method. 
Curvature estimates for starshaped hypersurfaces with given $k$-symmetric curvature have 
been established for various ambient Riemannian manifolds. For hypersurfaces in the sphere, the lower order and the curvature estimate are given in \cite{Oliker1} by M. Barbosa, L. Herbert and V. Oliker. These were used to prove the existence result by Y. Li and V. Oliker in \cite{YYLi:elli} via a degree theory argument. The curvature estimate and the existence result for prescribed curvature hypersurfaces in the hyperbolic space was proved by Q. Jin and Y. Li in \cite{YYLi:hyp} using similar arguments of W. Sheng, J. Urbas and X. Wang in \cite{sheng2004}. 

\noindent
For spacelike hypersurfaces in Lorentz manifolds, less is known. The case of prescribed mean curvature was studied by C. Gerhardt \cite{Gerhardt1983, Gerhardt1997} and R. Bartnik and L. Simon \cite{bartnik1982}.  C. Gerhardt \cite{CG:Lorentz} also studied more general curvature functionals, but was forced to only consider a comparatively restrictive class of curvature functionals due to a problematic sign on a term in the curvature estimates (in particular all of these functionals are zero on the boundary of the positive cone, 
 $\Gamma^+$). Urbas \cite{Urbas:cespacelike} obtained curvature estimates for prescribed symmetric curvature under the assumption that the mean curvature was in $L^p$ for sufficiently large $p$. Y. Huang \cite{Huang:01} noted that over compact domains Minkowski space, if the prescription function is additionally required to depend on the tilt and satisfies certain structural assumptions, then the problematic term in the curvature estimates may be cancelled  and curvature estimates may be obtained. D. Ballesteros-Ch\'avez \cite{Ballesteros2019} extended this result to compact domains in de Sitter space.

 \noindent
We will prove the following:
\begin{thm}\label{FullTheorem}
Suppose that $\psi:S^{n+1}_1\times\bb{R}_+\ra\bb{R}$ is a smooth positive function satisfying the structural conditions \ref{BarriersE}--\ref{convex} and let $1\leq k\leq n$. Then there exists a smooth embedded $k$-admissable spacelike hypersurface $\Sigma\subset S^{n+1}_1$ satisfying (\ref{ActualEquation}).
\end{thm}
\noindent
In Section \ref{Preliminaries} we will collect all required definitions and some preliminary calculations. In Sections \ref{sec:c0-estimate} and \ref{sectiontilt} we will prove the required $C^0$ and $C^1$ estimates respectively. In Section \ref{sec:curvests} we extend earlier results \cite{Ballesteros2019,Huang:01} to give the curvature estimates (and therefore $C^2$ estimates). In Section \ref{sec:exists} we prove existence of a solution via regularity result of Evans and Krylov  \cite{Evans, Krylov} and the degree theory of Y. Li \cite{YYLi:1989}.  

\noindent
We now state our structural assumptions on the prescription $\psi$:
\begin{assump*}
We impose the following structural assumptions on $\psi:S^{n+1}_1\times[1.\infty)\rightarrow \mathbb{R}$: 
\begin{enumerate}[label=\textbf{\Alph*)}]
\item \label{BarriersE}(Barrier conditions) There exist constants $0<R_1<R_2<\infty$ such that
  \begin{flalign*}
       \tanh(r)>\psi(Y(r, \xi), \cosh(r))&\qquad\quad \mbox{for all } \xi\in\mathbb{S}^n, r< R_1,\\
      \tanh(r)<\psi(Y(r, \xi), \cosh(r))&\quad\qquad\mbox{for all } \xi\in\mathbb{S}^n, r > R_2.
      \end{flalign*}
      where $Y$ is as in equation (\ref{eq:37}).
\item (Differential inequality)  For all $x\in S^{n+1}_1$ and $\tau\in [1,\infty)$,\label{differentialineq}
\[\psi_\tau(x,\tau) \tau\geq\psi(x,\tau)\]
\item (Asymptotics) $\frac{\psi(x,\tau)}{\tau} \rightarrow \infty$ as $\tau \rightarrow \infty$ for all $x\in S^{n+1}_1$.\label{asymptotics}
\item ($C^1$ bound) Taking coordinates $(x^1, \ldots, x^{n+1})=(r, \xi_1, \ldots, \xi_n)$ on $S^{n+1}_1$ as in (\ref{eq:37}) we have that there exists a uniform $C>0$ such that for all $x\in S^{n+1}_1$ and $\tau\in\bb{R}$,\label{estimates}
\[\left|\pard{\psi(x,\tau)}{x^i}\right|\leq C \psi(x,\tau)\ .\]
\item (Convexity in $\tau$) $\psi_{\tau\tau}(x,\tau)\geq 0$ for all $x\in S^{n+1}_1$ and $\tau\in [1,\infty)$.\label{convex}
\end{enumerate}
\end{assump*}
\begin{rem}
Assumption \ref{BarriersE} is simply to ensure the existence of suitable barriers. The relevance of the rate $\tanh(r)$ is that this is curvature of the natural foliation of totally umbillic hypersurfaces in de Sitter space.
\end{rem}
\begin{rem} We note that condition \ref{differentialineq} already implies that $\psi(x,\tau)\geq\tau\psi(x,1)$, so \ref{asymptotics} may be considered as the smallest possible increase in growth on top of this assumption.
\end{rem} 
\begin{rem} Condition \ref{estimates} above is used to estimate the space derivative of $\psi$ with respect to a multiple of $\psi$. This is a necessary condition in our tilt estimates as the derivative may be vastly larger with respect to $\tau$, for example if in local coordinates $\psi(x,\tau)=\tau^2+x^1e^\tau$ then at $x^1=0$, $\pard{\psi(x,\tau)}{x^1}$ cannot be estimated by $\psi$.
\end{rem}
\begin{rem}
There are an abundance of functions $\psi$ which satisfy the structure conditions \ref{BarriersE}--\ref{convex} (see also Lemma \ref{BarriersExist}). Our model function is $\psi(x,\tau) = \Psi(x)\tau^p$ where $p>1$ and $\Psi$ is a smooth bounded function satisfying the conditions of Lemma \ref{BarriersExist}. 
\end{rem}

\section{Preliminaries}
\label{Preliminaries}

\label{sec:geom-form-hypers}

\subsection{Subspace geometry in {L}orentzian manifolds}
To avoid confusion with signs, we now collect some geometric formulae for hypersurfaces in Lorentzian manifolds.
\noindent
Let $\{\partial_1,...,\partial_n,N\}$ be a basis for a Lorentzian
manifold $(\ov{M},\bar{g})$ and $M$ a Lorentzian  (not necessarily spacelike) hypersurface with induced metric $g$ such that
$\{\partial_i\}$ span $TM$, and let $N$ be the unit normal field to $M$
and put $\epsilon = \bar{g}(N,N)$. When the induced metric is positive
definite, then we say that $M$ is a spacelike hypersurface, then $g$
can be represented by the matrix
$g_{ij} = g(\partial_i,\partial_j)$ with inverse denoted by $g^{ij}$.
\noindent
The \textit{Gauss formula} for $X,Y\in T\Sigma$ reads
\begin{equation*}
  \label{eq:4}
  D_{X}Y = \nabla_{X}Y +\epsilon\, h(X,Y)N,
\end{equation*}
here $D$ is the connection on $\ov{M}$, $\nabla$ is the induced connection on $M$ and the \textit{second fundamental form} $h$ is 
the normal projection of $D$. In coordinate basis we write
\begin{equation*}
  \label{eq:5}
h_{ij}= h(\partial_i,\partial_j).
\end{equation*}
The \textit{shape operator} obtained by raising an index with the inverse of the metric
\begin{equation*}
  \label{eq:6}
  h^i_j = g^{ik}h_{kj}.
\end{equation*}
The \textit{principal curvatures} of the hypersurface $\Sigma$ are the
eigenvalues of the symmetric matrix $(h^i_j)$. The tangential projection of the
covariant derivative of the normal vector field $N$ on $\Sigma$, $\nabla_j N = (D_{\partial_j}N)^{\top}$,
is related to the second fundamental form by the \textit{Weingarten equation}
      \begin{equation}
        \label{eq:7}
        \nabla_{j} N = -h^i_j\partial_i = -g^{ik}h_{kj}\partial_i.
      \end{equation}
      The \emph{curvature tensor} is defined for $X,Y,Z \in T\Sigma$ as
      \begin{equation*}
        \label{eq:8}
        R(X,Y)Z = \nabla_{Y}\nabla_{X}Z - \nabla_{X}\nabla_{Y}Z + \nabla_{[X,Y]}Z.
      \end{equation*}
      Contracting with the metric
      \begin{equation*}
        \label{eq:11}
        R_{ijkl} = g\left(R(\partial_i,\partial_j)\partial_k,\partial_l\right) = g_{lm}R^{m}_{ijk}.
        \end{equation*}
  To relate curvature and second fundamental form, we have the \emph{Codazzi equations}
 \begin{equation}
  \label{GeneralCodazzi}
  \ip{\ov{R}_{ijk}}{N} = \nabla_j h_{ik}-\nabla_i h_{jk}.
\end{equation}     
and the \textit{the Gauss equation},
\noindent
\begin{equation}
  \label{eq:15}
  \ov{R}_{ijkl} = R_{ijkl} - \epsilon \left(h_{ik}h_{jl} - h_{il}h_{jk} \right).
\end{equation}

\noindent
Note that if $T$ is a symmetric tensor, then the following \textit{Ricci identity} holds
      \begin{equation}
        \label{eq:17}
{\nabla_{k}}\nabla_{l}T_{ij} - {\nabla_{l}}\nabla_{k}T_{ij} =  R_{kljr}T_{ir} + R_{klir}T_{rj}.
      \end{equation}
      
\subsection{The geometry of de {S}itter space}
We now consider manifolds $\Sigma\subset S_1^{n+1}\subset\mathbb{R}^{n+2}_1$ where 
\begin{itemize}
    \item  $\mathbb{R}^{n+2}_{1}=(\mathbb{R}^{n+2},\bar{g})$ is \emph{Minkowski space} with metric \[\bar{\bar{g}} = -dx_1^2 + dx_2^2 + \cdots + dx_{n+2}^2\] and
covariant derivative $\bar D$. 
\item $S^{n+1}_1$ is \emph{de Sitter space}, defined by
\[{S}^{n+1}_{1} = \left\{ x\in\mathbb{R}^{n+2}_{1} : -x_1^2 + x_2^2 + \cdots + x_{n+2}^2 = 1\right\}\] with the
induced Lorentzian metric $\ov{g}$, covariant derivative $D$, unit normal $N$  and second fundamental form $h$.
\item $\Sigma \subset S^{n+1}_1$ is a embedded spacelike hypersurface of $S^{n+1}_1$ with induced Riemannian metric $g$, covariant derivative $\n$, unit normal $\nu$ and second fundamental form $A$.
\end{itemize}
Let $\bb{S}^n$ be the standard round sphere. Then de Sitter space may be parametrised by $Y:\bb{S}^n\times \bb{R}\ra S^{n+1}_1$ given by
\begin{equation}
             \label{eq:37}
             Y(r,\xi) = \sinh(r)E_1 + \cosh(r) \xi
           \end{equation}
and in these coordinates, 
the induced metric is
\begin{equation*}
  \label{eq:18}
  \ov{g} = -dr^2 + \cosh^2(r)\sigma,
\end{equation*}
where $\sigma$ is the round metric on $\mathbb{S}^n$. 
\noindent
We note that as $|Y|^2=1$ we have that $Y=N$ is a unit normal to $S^{n+1}_1$ and (as $0=\ip{Y_{\alpha\beta}}{Y}+\ip{Y_\alpha}{Y_\beta}$)
\begin{equation}h_{\alpha\beta}=-\ov{g}_{\alpha\beta} \ ,\label{hisg}\end{equation}
(where $1\leq \alpha, \beta\leq n+1$) and so using the Gauss equation we have
\[R^{S^{n+1}_1}_{\alpha \beta\gamma\delta} = \ov{g}_{\alpha\gamma}\ov{g}_{\beta\delta}-\ov{g}_{\alpha\delta}\ov{g}_{\beta\gamma}\ ,\]
where we used that Minkowski space is flat. This implies that on $\Sigma$ (which has a timelike normal) we have that
\begin{equation}
\n_ih_{jk}=\n_jh_{ik}\label{CM}
\end{equation}
and 
\begin{equation}R_{ijkl}=A_{il}A_{jk}-A_{ik}A_{jl}+{g}_{ik}{g}_{jl}-{g}_{il}{g}_{jk}\ .\label{GE}\end{equation}
We define the \emph{tilt function} on $\Sigma$ to be the function
\begin{equation}\tau = \langle \nu , E_1 \rangle ,\label{Tiltdefinition}\end{equation}
where $\nu$ is a unit normal to $\Sigma$ which has been chosen so that $\tau$ is positive.
\noindent
We now represent $\Sigma$ as a graph, that is we take $u:\bb{S}^n\ra\bb{R}$ so that $\Sigma$ is parametrised by $X:\bb{S}^n\ra S^{n+1}_1$ given by $X(\xi)=Y(u(\xi),\xi)$. We will use $\tilde{\nabla}$ to denote the standard covariant derivative for the metric $\sigma$ on $\bb{S}^n$, and our indices $\partial_{i},\partial_j,...$, etc., take values from $1$ to $n$, 
except for the vector field $\partial_{r}$ which will be considered separately. 
The tangent space of the hypersurface at a point $Y\in\Sigma$
is spanned by the tangent vectors
$Y_j = u_{j}\partial_{r} + \partial_{j}$, the covariant derivative
 $\nabla$ corresponding to the induced metric 
on $\Sigma$  is given by
\begin{equation*}
  \label{eq:19}
  g_{ij} = -u_{i} u_{j} + \cosh^2(u)\sigma_{ij} .
\end{equation*}
We write 
\[\tilde{\tau}=\frac{\cosh^2(u)}{\sqrt{\cosh^2(u)-|\tilde{\n} u|^2}}\ ,\]
where $\tilde{\nabla}u = \sigma^{ij}u_j\partial_i$ and $|\tilde{\nabla}u| := \sigma^{ij}u_iu_j$ (we will see shortly that $\tilde{\tau}=\tau$). $\Sigma$ is spacelike at a point if $g_{ij}$ is invertable, which is equivalent to $\tilde{\tau}$ being finite at that point. 
Since $\Sigma$ is spacelike, we calculate the inverse of $g_{ij}$ to be
  \begin{equation}
    \label{eq:ginverse}
    \begin{split}
      g^{ij}
    &=\cosh^{-2}(u) \sigma^{ij} +
    \frac{ \sigma^{il}u_{l} \sigma^{jm}u_{m}}{\cosh^4(u) - \cosh^{2}(u)|\tilde{\nabla} u|^2}=\cosh^{-2}(u)\left[\sigma^{ij} +\tilde{\tau}^2 \frac{\sigma^{il}u_{l} \sigma^{jm}u_{m}}{\cosh^4(u)}\right]
    \end{split}
  \end{equation}
 A unit normal vector
  to $\Sigma$ at the point $Y$ can be obtained by solving the equation $\ov{g}(Y_{\alpha},\hat{n})=0$, and then we get
\begin{equation*}
  \label{eq:21}
  \nu =- \frac{ \cosh^2(u)\partial_{r}+ \tilde{\nabla} u }{
    \sqrt{\cosh^4(u) - \cosh^2(u)|\tilde{\nabla} u|^2}}=\cosh^{-3}(u)\tilde{\tau}(\cosh^2(u)\partial_r+\tilde \n u).
\end{equation*}
We note that 
           \begin{equation*}
             \label{eq:38}
             \cosh(u)\partial_r =
             E_1 + \sinh(u)Y.
           \end{equation*}
and so we see that
\[\tau=\ip{E_1}{\nu}=\tilde{\tau}\ ,\]
as claimed.
\noindent
The second fundamental form is the projection of the second derivatives of the
parametrisation $D_{Y_{\alpha}}Y_{\beta}$ on the normal direction. Writing $\tilde{\Gamma}$ for the Christoffel symbols of the metric $\sigma$, we have
  \begin{equation*}
    \label{eq:23}
    D_{\partial_r}\partial_r = 0; \quad
    D_{\partial_r}{\partial_j} = \tanh(r)\partial_j;\quad
    D_{\partial_i}{\partial_j}= \cosh(r)\sinh(r)\sigma_{ij}\partial_{r}+ \tilde\Gamma_{ij}^k \partial_k ,
  \end{equation*}
and using these identities we compute
  \begin{equation*}
    \label{eq:24}
    \begin{split}
    D_{Y_i}Y_{j}& = D_{u_i\partial_{r}+\partial_i}\left(u_j\partial_{r}+\partial_j\right)\\
    &=
    u_ju_jD_{\partial_{r}}\partial_{r} +
    u_iD_{\partial_{r}}\partial_{j} +
    u_{ij} \partial_{r} +
    u_{j}D_{\partial_{i}}\partial_{r} +
    D_{\partial_{i}}\partial_{j}.
    \end{split}
  \end{equation*}
It follows that  $ A_{ij}=  \ov{g}(D_{Y_{i}}Y_{j}, \nu)$ is given explicitly by
 \begin{equation}
  \label{eq:Secondff}
  A_{ij} = \cosh^{-1}(u)\tau
  \left(\tilde{\nabla}^2_{ij}u - 2\tanh(u)u_{i}u_{j} + \sinh(u)\cosh(u)\sigma_{ij}\right).
\end{equation}

\noindent
Finally, we define a notion of partial derivatives for $\psi$ on $S^{n+1}_1$. Suppose we have a function $f:S^{n+1}_1\times \bb{R}\ra \bb{R}$, then we define the partial derivative on $S^{n+1}_1$ and $\Sigma$ by
\begin{align*}
D^x f =\pard{\psi}{x^\alpha}\ov{g}^{\alpha\beta}\partial_{\beta}, \qquad
\n^x \psi = (D^x f)^\top\ . 
\end{align*}
Similarly, we may define second derivatives of $f$ in the usual tensorial way.

\subsection{Curvature functionals and admissability}
Throughout this section we fix $1\leq k\leq n$. As described above we will consider functions
\[F[A] := H_k^\frac 1 k(\lambda[A]),\]
where $\lambda[A]=(\lambda_1,\ldots, \lambda_n)$ are the eigenvalues of the symmetric matrix $A$ and we will define $f:\bb{R}^n\ra\bb{R}$ by $f=H_k^\frac{1}{k}$. We define the \emph{admissable cone} of $H_k$ to be
$\Gamma_{k}$ which is defined to be the connected component of $H_k^{-1}(\{x\in\bb{R}|x>0\})$ which contains the positive cone
\[\Gamma^{+}=\{\lambda\in\mathbb{R}^n\,|\, \lambda_i>0, \forall i = 1,2,\dots,n\}.\]
We have that for all $\lambda\in\Gamma_k$, $f_{\lambda_i}(\lambda) >0$ and $f$ is concave in $\Gamma_k$. 
Since $f \in C^2(\Gamma_k)\cap C^{0}(\ov{\Gamma_k})$ it follows that $F[A]$ is elliptic and concave if the eigenvalues of $A$ lie in $\Gamma_k$.
\noindent
A spacelike hypersurface $\Sigma\subset S^{n+1}_1$ will be called \emph{$k$-admissable} if for all $p\in\Sigma$ the eigenvalues of the shape operator $A_i^j=A_{ik}g^{kj}$ are in $\Gamma_k$.
\noindent
A function $u\in C^2(\bb{S}^n)$ will be called \emph{$k$-admissable} if the graph of $u$ is a spacelike admissable hypersurface. We note that this implies that $u$ is positive, as if there is a negative minimum of $u$ at $p\in\Sigma$, the shape operator is negative definite so $\lambda[A]\notin \Gamma_k$.
\noindent
As is standard, and we will write
\[F^{ij} := \pard{F}{A_{ij}}, \qquad\qquad F^{ij,kl} := \frac{\partial^2F}{\partial A_{ij}\partial A_{kl}}\ .\]

\section{A priori $C^0$ estimate}
\label{sec:c0-estimate}
\noindent
Considering $\Sigma$ graphically, any solution $u\in C^2\left(\mathbb{S}^n\right)$ of \eqref{ActualEquation}, the barrier conditions \eqref{barrier} will ensure that $R_1 \leq u(\xi) \leq R_2$ for all $\xi\in\mathbb{S}^n$. The proof follows maximum principle arguments similar to \cite[Lemma 3.1]{Oliker1}.
\begin{lem}\label{BarrierConditions}
  Let $1\leq k\leq n$ and let $\psi: (0,\infty)\times \mathbb{S}^n\times[1, \infty)\to \mathbb{R} $ is a continuous positive function such that there exist constants $0<R_1<R_2<\infty$ such that
  \begin{equation}
    \label{barrier}
    \begin{split}
       \tanh(r)>\psi(r, \xi, \cosh(r)),\quad \mbox{for } \xi\in\mathbb{S}^n, r< R_1,\\
      \tanh(r)<\psi(r, \xi, \cosh(r)),\quad\mbox{for } \xi\in\mathbb{S}^n, r > R_2.
      \end{split}
    \end{equation}
    Then $u\in C^2\left(\mathbb{S}^n\right)$ is a solution of \eqref{ActualEquation} then for all $\xi\in\mathbb{S}^n$, $u$ satisfies
    \begin{equation}
      \label{eq:25}
      R_1 \leq u(\xi) \leq R_2.
    \end{equation}
\end{lem}
  \begin{proof}
    Suppose there is a point $\xi_0\in\mathbb{S}^n$ where the maximum of $u$ is attained, say $R_2<r_0 = u(\xi_0)$. Then at the maximum we have $\tilde{\nabla} u \dot{=} 0$ and $\tilde{\nabla}^2 u \dot{\leq} 0$. Substituting these into the equations of the previous section Then the inverse of the  metric, the tilt and the second fundamental form at $\xi_0$ are respectively
      \begin{equation*}
        \label{eq:26}
        g^{ij} \dot{=}\frac{1}{\cosh^2(r_0)}\sigma^{ij}, \qquad \tau \dot{ =}  \cosh u(u)=\cosh(r_0), \qquad
        A_{ij}\dot{\leq} \sinh(r_0)\cosh(r_0)\sigma_{ij}.
      \end{equation*}
    At $\xi_0$, the shape operator therefore satisfies
      \begin{equation}
        \label{eq:28}
        A^i_j \leq \tanh(r_0) \delta^i_j,
      \end{equation}
      and so $\lambda_i\leq\tanh (r_0)$ for $1\leq i \leq n$. Substituting into \eqref{eq:44}, 
      \begin{equation}
        \label{eq:29}
       \psi(r_0,\xi_0, \cosh(r_0)) \dot{=} F[A^i_j] \leq \tanh(r_0) < \psi(r_0,\xi_0, \cosh(r_0)),
      \end{equation}
     which is a contradiction, and so  $u \leq R_2$. An analogous argument at the minimum completes the proof.
  \end{proof}
\noindent
Observe that, we may impose a few fairly mild assumptions on $\psi_r$ to ensure that barriers exist.
\begin{lem}
 Let $\psi: [0,\infty)\times \mathbb{S}^n\times[1, \infty)\to \mathbb{R}$ be a uniformly bounded in $C^2$ which is positive for $r>0$. We give conditions for upper and lower barriers:\label{BarriersExist}
\begin{description}[leftmargin=*]
 \item[Lower Barriers]
 If for all $\xi\in \mathbb{S}^n$ 
\[\psi(0,\xi, 1)=0 \qquad\text{ and }\qquad \psi_r(0,\xi,1)<1\] 
 then a lower barrier exists, that is there exists an $0<R_1$ such that for all $\xi\in\mathbb{S}^n$ and $r<R_1$, 
 \[\tanh(r)>\psi(r, \xi, \cosh(r))\ .\]
\item[Upper Barriers]  Suppose $\psi$ satisfies the structural condition $\psi_\tau\tau>\psi$ and let $\gamma(r)$ be any smooth monotonic function s.t. $\gamma(r)\rightarrow\infty$ as $r\rightarrow\infty$. Suppose additionally that for all $\xi\in \mathbb{S}^n$ and $r>\tilde{R}$,
\[ \psi_r(r,\xi,\cosh(r))>\left[\gamma'(r)-\tanh(r)\right]\psi\ ,\] 
then an upper barrier exists, that is there exists an $0<R_2$ such that for all $\xi\in\mathbb{S}^n$ and $r>R_2$, 
 \[\tanh(r)<\psi(r, \xi, \cosh(r))\ .\]
\end{description} 
\end{lem}
 \begin{proof}
 We fix $\xi$ and consider the function
 \[\Psi(r,\xi) = \frac{\psi(r, \xi, \cosh (r))}{\tanh(r)}\ .\]
 \noindent
 Finding a lower barrier is equivalent to showing that for all $\xi\in\mathbb{S}^n$ there exists a $R_1>0$ such that for all $0<r<R_1$, $\Psi(r,\xi)<1$. By the assumptions, $\Psi$ is continuous up to $r=0$, and so by compactness, there exists a $\delta>0$ such that for all $\xi\in \mathbb{S}^n$, $\Psi(0,\xi)<1-\delta$. By continuity and compactness there exists an $R_1$ such that $\Psi(r,\xi)<1$ for all $0<r<R_1$ as claimed.
 \noindent
  Finding an upper barrier is equivalent to showing that for all $\xi\in\mathbb{S}^n$ there exists a $R_2>0$ such that for all $r>R_2$, $\Psi(r,\xi)>1$. We calculate that
  \begin{flalign*}
  \frac{d}{dr}\Psi(r,\xi) &= -\frac{1}{\sinh^2(r)}\psi(r,\xi,\cosh(r)) +\coth(r)(\psi_r+\psi_\tau \sinh(r))\\
  &\geq\left[1-\frac{1}{\sinh^2(r)}\right]\psi +\coth(r)\psi_r\\
  &\geq \frac{\gamma'(r)}{\tanh(r)}\Psi
  \end{flalign*}
 for $r>\hat{R}(\tilde{R})$ sufficiently large depending on $\tilde{R}$. Since $\tanh(r)<1$, we seee that \[\frac{d}{dr}\Psi(r)>\gamma'(r)\Psi\ ,\]
 which implies $\Psi(\xi,r)>\Psi(\hat{R}, \xi)e^{\gamma(r)}$. Since $\mathbb{S}^n$ is compact, and $\gamma(r)\rightarrow\infty$ as $r \rightarrow\infty$, this implies that the claimed upper barrier conditions are met.
 \end{proof}
 
 \begin{rem}
 A function $\psi$ satisfying lower barrier conditions in Lemma \ref{BarriersExist} allows the solution to \eqref{eq:44} given by $u\equiv0$. As noted earlier, this solution is inadmissable, as we require strictly positive $u$ for the shape operator to be in the admissable cone everywhere.
 \end{rem}

\section{Tilt estimate}
\label{sectiontilt}
\noindent
We now demonstrate a strict spacelikeness estimate by estimating the tilt function $\tau$. The height function is defined by
           \begin{equation}
             \eta := -\langle Y, E_1\rangle, \label{etadef}
           \end{equation}
and we note that in terms of the graph function, $\eta=\sinh(u)$.
We now demonstrate the following identities.
           \begin{prop}
             \label{prop:1}
             The tilt and height functions satisfy the following identities:
             \begin{enumerate}
             \item\label{p:1}{$\nabla^2_{ij}\eta = \tau A_{ij} + \eta g_{ij}$.}
             \item\label{p:2}{$\nabla_j \tau = A_{j}^i\nabla_i\eta$,}
             \item\label{p:3}{$\nabla_j\nabla_i \tau =  \nabla_{k}A_{ij}\nabla^k\eta + \tau  A^2_{ij} +  A_{ij} \eta$,}
             \end{enumerate}
             where $A^2_{ij}:=A_{kj}A^k_i$.
             \end{prop}
             
    \begin{proof}
           Using (\ref{hisg}) and the Gauss formula we have that
           \begin{equation*}
             \begin{split}
             \label{eq:42}
             \nabla^2_{ij}\eta & = Y_{j}(Y_{i}\eta) - (\nabla_{Y_{i}}Y_{j}) \eta\\
             & = Y_{j}(Y_{i}\eta) - (\bar{D}_{Y_{i}}Y_{j}-h_{ij}N + A_{ij}\hat{n}) \eta\\
             & = \bar{D}^2_{Y_{j}Y_{i}}\eta + (h_{ij}N - A_{ij}\nu) \eta\\
             & =\tau A_{ij} + \eta g_{ij},
             \end{split}
           \end{equation*}
where we used that $X(\eta) = -\ip{X}{E_1}$ and $\bar{D}^2\eta=0$.
 \noindent
               Using the Weingarten equation (\ref{eq:7}) we obtain
           \begin{multline*}
             \label{eq:41}
             \nabla_j \tau = \langle \nabla_j \nu, E_1\rangle
             = -g^{ik}A_{kj}\langle Y_i, E_1\rangle
             = -g^{ik}A_{kj}\nabla_i\langle Y, E_1\rangle
             = g^{ik}A_{kj}\nabla_i\eta.
           \end{multline*}
           and from this we have
           \begin{equation*}
             \begin{split}
             \label{eq:43}
             \nabla^2_{ij}\tau &= \nabla_{j}(g^{mn}A_{ni}\nabla_{m}\eta)\\
             & = g^{mn}\nabla_{j}A_{ni}\nabla_{m}\eta + g^{mn}A_{ni}\nabla_{mj}\eta\\
             & = g^{mn}\nabla_{n}A_{ij}\nabla_{m}\eta + \tau  A_{mj}g^{mn}A_{ni} +  A_{ij} \eta
           \end{split}
         \end{equation*}
         where we used (\ref{CM}) on the third line.
               \end{proof}

\begin{prop}\label{TiltEstimate}
Suppose that $\psi:S^{n+1}\times\bb{R}\ra\bb{R}$ is smooth and positive and satisfies assumptions \ref{asymptotics} and \ref{estimates} above. Suppose that $u\in C^3(\bb{S}^n)$ satisfies (\ref{ActualEquation}) so that there exist $R_1, R_2 >0$ such that,
\[0<R_1<u(\xi)<R_2\]
for all $\xi\in\mathbb{S}^n$. Then there exists a constant $C_\tau$, depending only on $n$, $k$, $R_1$,$R_2$ and $\psi$ such that
\[\tau<C_\tau\ .\]
\end{prop}

\begin{proof}
We have that
\begin{flalign*}F^{ij}\nabla_j\nabla_i \tau &=  \nabla^k F\nabla_{k}\eta + \tau  F^{ij}A_{mj}A^m_{i} +  F \eta\\
&=\ip{\nabla^x \psi}{\nabla \eta} + \psi_\tau \ip{\nabla \tau}{\nabla \eta}+\tau F^{ij}A^2_{ij} +\psi \eta
\end{flalign*}
Due to structural assumption \ref{estimates} on $\psi$, we have that
\[|\ip{\nabla^x \psi}{\nabla\eta}|\leq C\tau^2\psi .\]
Furthermore, following \cite[equation 3.8]{Urbas:cespacelike} and using Newton's inequalities \cite{NewtonInequalities} (which are valid for $S_k$ \emph{outside} $\Gamma_k$), 
\begin{equation}F^{ij}A^2_{ij} \geq  H_k^\frac 1 k H_1\geq H_k^\frac{2}{k}=\psi^2\ .\label{UrbasMagicInequality}\end{equation}
Substituting this into the above equation we see that at a maximum of $\tau$, as $\n \tau=0$ and $\n_i\n_j\tau\leq 0$, we have that 
\begin{flalign*}
0&\geq - C\tau^2\psi +\psi \eta +\tau\psi^2
\end{flalign*}
which implies, by using the $C^0$ estimates that
\[\psi \leq C \tau\ .\]
The $C^0$ estimates imply that the solution stays in a compact region of $S^{n+1}_1$ and so, due to structural assumption \ref{asymptotics} on $\psi$, there exists a uniform $\tau_0$ such that for all $x$ in this region and for all $\tau>\tau_0$, $\tau<\frac{\psi}{2C}$. This yields a contradiction to $\tau>\tau_0$ and proves the Lemma.
\end{proof}

\section{A priori $C^2$ estimate}
\label{sec:curvests}
\noindent
From \cite{Ballesteros2019} we have the following curvature estimates over domains $\Omega \subset S^n$.
               \begin{thm}[Ballesteros-Ch\'avez 2019]
                 \label{thm:01}
                 Let $\Omega\subset\mathbb{S}^n$ be a domain in the round sphere, and let $u\in C^{4}(\Omega) \cap C^2(\bar{\Omega})$ an admissible solution of
                 the boundary value problem
                 \begin{equation*}
                   \label{eq:44}
                   \left\{
                     \begin{array}{rcll}
                      F(A) =  H_{k}^{\frac{1}{k}}(\lambda(A)) &=& \psi(Y,\tau) & \mbox{in}\quad \Omega\\
                       u  &=& \varphi & \mbox{on}\quad \partial\Omega
                     \end{array}
                     \right.,
                   \end{equation*}
                   where $A$ is the second fundamental form of a spacelike surface $\Sigma$ in de Sitter space given by
                   (\ref{eq:Secondff}) and $\psi$ is a smooth positive function satisfying \ref{differentialineq} and \ref{convex}. Assume additionally that there exists a $R_1,R_2, C_\tau>0$ such that  
     \[R_1<u(\xi)<R_2, \qquad \tau(\xi)<C_\tau .\]
                   Then
                   \begin{equation*}
                     \label{eq:46}
                     \sup_{\Omega}|A| \leq C,
                   \end{equation*}
                   where $C$ depends on $n$, $|\varphi|_{C^1(\bar{\Omega})}$,$R_1$, $R_2$, $C_\tau$, $|\psi|_{C^2([R_1,R_2]\times\Omega\times[1,C_{\tau}])}$ and $\sup_{\partial\Omega}|A|$.
                 \end{thm}
\noindent  
We will now extend this result to all of $\bb{S}^n$, or equivalently obtain estimates on all of $\Sigma$. 
\begin{thm}\label{CurvatureEstimate}
Suppose that $\psi$ is a smooth positive function which satisfies \ref{differentialineq} and \ref{convex}. Suppose $u\in C^4(\bb{S}^n)$ is a solution of (\ref{ActualEquation}) such that there exist constants $0<R_1<R_2$ and $C_\tau>0$ such that for all $\xi\in \bb{S}^2$,
\[R_1<u(\xi)<R_2, \qquad \tau(\xi)<C_\tau .\]
Then there exists a constant $C_A=C_A(k,n,R_1,R_2,C_\tau,|\psi|_{C^2([R_1,R_2]\times\bb{S}^n\times[1,C_\tau))})$ so that
\[|A|<C_A\ .\]
\end{thm}
\begin{proof}
Suppose first that $k=1$. In this case (\ref{eq:ginverse}) and (\ref{eq:Secondff}) imply that if we write (\ref{ActualEquation}) in terms of the graph function $u$, we obtain a quasilinear equation which is uniformly elliptic if we have a uniform bound on $\tau$. Therefore if we have the assumed bounds, we see that this equation is uniform ellipticity with uniform $C^1$ estimates on the solution. We may therefore apply De Giorgi--Nash--Moser estimates and Schauder estimates to imply uniform $C^2$ estimates. In this case, the theorem is therefore proven.

\noindent
Suppose now that $k\geq 2$. We begin by proving a Simon's-type identity for the second fundamental form. At an arbitrary point $p\in\Sigma$, we choose coordinates so that $A_{ij}$ is diagonal and $g_{ij}=\delta_{ij}$. In these coordinates $F^{ij}$ is also diagonal.
\noindent
The Codazzi equation (\ref{CM}) and the Ricci identity (\ref{eq:17}) imply that 
\begin{equation*}
             \begin{split}
               \nabla_{i}\nabla_{j}A_{kk} =& \nabla_{i}\nabla_{k}A_{kj} \\
               =&  \nabla_{k}\nabla_{i} A_{kj} + R_{ikkr}A^r_{j} + R_{ikjr}A_{k}^r\\
               =& \nabla_{k}\nabla_{k} A_{ij} + R_{ikkr}A_{j}^r + R_{ikjr}A_{k}^r.
               \end{split}
           \end{equation*}
where we only sum over indices where one is raised and one is lowered. Choosing $i=j$ in the above, and applying the Gauss equation (\ref{GE}) we obtain
\begin{align*}
\nabla_{j}\nabla_{j}A_{kk} =& \nabla_{k}\nabla_{k} A_{jj} + R_{jkkr}A_{j}^r + R_{jkjr}A_{k}^r\\
&=\nabla_{k}\nabla_{k} A_{jj} + (A_{jr}A_{kk} - A_{jk}A_{kr}+g_{jk}g_{kr}-g_{jr}g_{kk})A_{j}^r\\
&\qquad +(A_{jr}A_{kj}-A_{jj}A_{kr}+g_{jj}g_{kr}-g_{jr}g_{kj})A_{k}^r\\
&=\nabla_{k}\nabla_{k} A_{jj} + A^2_{jj}A_{kk} - A_{jk}A^2_{kj}+g_{jk}A_{jk}-A_{jj}g_{kk}\\
&\qquad +A^2_{jk}A_{kj}-A_{jj}A^2_{kk}+g_{jj}A_{kk}-A_{kj}g_{kj}\\
&=\nabla_{k}\nabla_{k} A_{jj} + A^2_{jj}A_{kk} -A_{jj}g_{kk}-A_{jj}A^2_{kk}+g_{jj}A_{kk}\ .
\end{align*}

\noindent
We therefore see that
\begin{flalign*}
F^{ij}\n_i\n_j A_{kk}&=\sum_{j=1}^nF^{jj}\n_j\n_j A_kk\\
&=F^{ij}\n_k\n_kA_{ij}+A_{kk}F^{ij}A^2_{ij}-F^{ij}A_{ij}A^2_{kk}+A_{kk}\text{tr}{F^{ij}}-F^{ij}A_{ij}\\
&=F^{ij}\n_k\n_kA_{ij}+A_{kk}F^{ij}A^2_{ij}-\psi A^2_{kk}+A_{kk}\text{tr}{F^{ij}}-\psi
\end{flalign*}
and so writing $H=nH_1 = \sum_{k=1}^n A_{kk}$ we have that
\begin{equation}
    F^{ij}\n_i\n_j H  = F^{ij}\Delta A_{ij}+H F^{ij}A^2_{ij} + H \text{tr} F^{ij} -\psi\left(n+|A|^2\right)
\end{equation}
where $\Delta$ is the Laplace--Beltrami operator. To estimate the first term on the right hand side we note that by differentiating (\ref{ActualEquation}) twice gives
\[\Delta \psi = F^{ij}\Delta A_{ij} + F^{ij,kl}\n_t A_{ij}\n^tA_{kl}\leq F^{ij}\Delta A_{ij}\]
where we used the well known concavity of $F$.
\noindent
We have that
\begin{flalign*}\n_l\n_k \psi &= \n^x_l\n^x_k \psi  +\n^x_k \psi_\tau \n_l \tau+\n^x_l \psi_\tau \n_k \tau +\psi_\tau \n_l\n_k \tau +\psi_{\tau \tau}\n_l \tau \n_k\tau\\
&=D_l^xD_k^x \psi - A_{lk}D^x_{\hat{n}}\psi+\n^x_l \psi_\tau \n_k \tau +\psi_\tau \n_l\n_k \tau+\psi_{\tau \tau}\n_l \tau \n_k\tau\ .
\end{flalign*}
Since $H^2 = 2S_2+|A|^2$, so if $\l[A]\in \Gamma_k$, $|A|<H$ and so
\begin{flalign*}\Delta \psi&=g^{kl}D_l^xD_k^x \psi - H D^x_{\hat{n}}\phi+2\ip{\n^x \psi_\tau}{\n \tau} +\psi_\tau \Delta \tau+\psi_{\tau \tau}|\n \tau|^2\\
&=\psi_\tau\left[\ip{\n \eta}{\n H}+\tau |A|^2+H\eta\right] +\psi_{\tau \tau}|\n \tau|^2- H D^x_{\hat{n}}\psi+2\ip{\n^x \psi_\tau}{\n \tau} +g^{kl}D_l^xD_k^x \psi\\
&\geq\psi_\tau\ip{\n \eta}{\n H}+\psi_\tau \tau |A|^2 - C_1 H-C_2\ ,
\end{flalign*}
where used structural assumption \ref{convex} of $\psi$ and we estimated using the bounds on $\tau$ and $u$ in compactness arguments to estimate derivatives of $\psi$. 
\noindent
Overall we have that
\begin{flalign*}
F^{ij}\n_i\n_j H & \geq H F^{ij}A^2_{ij} + H \text{tr} F^{ij}+(\psi_\tau\tau-\psi)|A|^2\\
&\qquad+\psi_\tau\ip{\n \eta}{\n H} -C_1 H -C_2-n\psi \ .
\end{flalign*}
Using (\ref{UrbasMagicInequality}) and the structural assumption \ref{differentialineq} on $\psi$ we have that
\begin{flalign*}
F^{ij}\n_i\n_j H &\geq \frac 1 n H^2\psi +H\text{tr}F^{ij}+\psi_\tau\ip{\n \eta}{\n H} -C_1 H -C_2-n\psi\ .
\end{flalign*}
Using the bounds on $u$, there exists a small constant $\delta>0$ such that $\psi>\delta$, and so at a maximum point of $H$,
\begin{flalign*}
0&\geq F^{ij}\n_i\n_j H\\
&\geq \frac \delta n H^2-C_1H-C_3 , 
\end{flalign*}
and so $H$ is bounded. This implies a uniform bound on $H$ and a uniform bound on $|A|$ now follows as $k\geq 2$.
\end{proof}

 \section{Proof of existence}
 \label{sec:exists}
\noindent
We now prove existence of solutions to \eqref{ActualEquation}, following the proof of V. Oliker and Y. Li \cite{YYLi:elli}. Throughout this section $\Sigma$ will be considered in graphical coordinates with graph function $u$. Consider for $0<\alpha<1$, the subset of functions in $C^{4,\alpha}(\mathbb{S}^n)$ which are $k$-admissible, denoted by $C^{4,\alpha}_{ad}(\mathbb{S}^n)$.
The idea is to consider a one parameter family of prescription functions $\psi_t$ where
\[\psi_t(\xi,u,\tau) := t\psi(\xi,u,\tau) +(1-t)\Psi(\xi,u,\tau) , \]
where $\Psi(\xi,u,\tau)$ is to be chosen shortly. We define  $\Phi:C^{4,\alpha}_{ad}(S^n)\times[0,1]\rightarrow C^{2,\alpha}$, by
\begin{equation}\Phi(u,t):=H_k^\frac{1}{k}(u_t)-\psi_t(\xi,u_t,\tau(u_t))\label{oneparamfamily}\end{equation} 
for all $t\in[0,1]$. We will apply degree theory to ensure that there exists at least one solution to $\Phi(u_t,t)=0$ for all $t\in[0,1]$. As in \cite{YYLi:elli}, to be able to apply the beautiful degree theory of Y. Li \cite{YYLi:1989}, we need to verify the following three steps:
\begin{description}
\item[Step 1a)] Show that there exists a unique solution to 
\[H_k^\frac 1 k(u_0) = \Psi(x,u_0, \tau(u_0))\ .\]
\item[Step 1b)] Show that at $u_0$ the the linearisation of $\Phi$ is invertible.
\item[Step 2] Define a suitable set of admissable functions and show that all solutions of (\ref{oneparamfamily}) stay in this set.
\item[Step 3] Verify that we may apply degree theory to ensure that the degree of $\Phi(\cdot,1)$ is not zero, and so a solution exists as claimed.
\end{description}

We choose $\Psi$ to be 
\[\Psi(\xi, u, \tau) = \tau^p u \tanh(u)\ .\]
for some $p>1$.
\begin{proof}[Proof of Step 1a)]
We may easily verify that a solution exists to $\Phi(u_0,0)=0$, by considering constant functions. The hypersurfaces corresponding to $u_0=\lambda$ are totally umbillic with principal curvatures $\tanh(\lambda)$, and so $H_k^\frac 1 k=\tanh(\lambda)$. We may see that on such a hypersurface, $\tau = \cosh(\lambda)$ and so if $\lambda$ satisfies $\lambda \cosh^p(\lambda) = 1$ then $u_0$ is a solution. Clearly such a $\lambda$ exists as, writing the continuous function $\varphi(x):=x\cosh^p(x)$, $\varphi(0)=0$ and $\varphi(1)>1$.

\noindent
Suppose that there exists another $u\in C^{4,\alpha}_a(S^n)$ satisfying $\Phi(u,0)=0$. Suppose furthermore that $\max u = u(\xi_0)=\lambda_0>\lambda$. As in the proof of the $C^0$ estimates we have that at $\xi_0$, 
\[\psi(\xi_0,u(\xi_0),\tau(\xi_0)) = \cosh^p(\lambda_0)\lambda_0\tanh(\lambda_0) =H_k^\frac 1 k|_{\xi_0}\leq\tanh(\xi_0)\]
which is a contradiction as $x\cosh^p(x)$ is a monotonically increasing function. Therefore $\max u \leq \lambda$. An identical argument implies $\min u\geq \lambda$, implying that $u(\xi)=u_0$. 
\end{proof}

\begin{proof}[Proof of step 1b)]
Considering $\tau$, $A_i^k$ as algebraic functions of $\xi$, $u$, $\tilde \n u$, $\tilde \n^2 u$, which we will write with the variables $x$, $r$, $p^i$, $z_{ij}$ respectively then $A_i^k=A_i^k(x,r,p,z)$ and $\tau = \tau(x,r,p)$. Then the linearisation of the above in direction $v$ is given by
\begin{flalign}
\nonumber 0&=\frac{d}{ds}\left(\left[H^\frac{1}{k}_k-\psi_t\right](u+sv) \right)\\
&=F^{i}_k\pard{A^k_i}{z_{ij}}\tilde\n_{ij}v + \left[F^{i}_k\pard{A^k_i}{p_{k}}-\psi_\tau\pard{\tau}{p_k}\right]\tilde\n_kv +\left[F^{i}_k\pard{A^k_i}{r}-\psi_r - \psi_\tau\pard{\tau}{r}\right]v\ , \label{genlin}
\end{flalign}
From equations (\ref{eq:ginverse}) and (\ref{eq:Secondff}) we have that
\begin{align*}
A_i^k &=\frac {\tau}{\cosh^3(u)}\left(\sigma^{ij} + \tau^2\frac{u_l\sigma^{li}u_m\sigma^{mj}}{\cosh^4(u)}\right)\left(\tilde\n_{ij} u - 2\tanh(u)u_iu_j+\sinh(u)\cosh(u)\sigma_{ij}\right)
\end{align*}
and so we see that 
\begin{align*}
\pard{A^k_i}{r}\big|_{(x,u,\tilde \n u, \tilde \n^2 u)}&=\frac{\pard{\tau}{r}}{\tau}A_i^k-3\tanh(u)A_i^k+\tau\frac{\cosh^2(u)+\sinh^2(u)}{\cosh^3(u)}\delta_i^k \\&\qquad\qquad+ u_lu_m\left[\text{ bounded terms }\right]\ ,
\end{align*}
and
\[\pard{\tau}{r}\big|_{(x,u,\tilde \n u, \tilde \n^2 u)}=\tanh(u)\tau +|\tilde \n u|^2\left[\text{ bounded terms }\right]\ .\]
\noindent
We are interested in (\ref{genlin}) when $t=0$, that is when $u=\l>0$ and $\psi = \Psi = \tau^pu\tanh(u)$. In this case, we have that
\[\tau = \cosh(u), \qquad A_i^k = \tanh(u) \delta_i^k,\qquad H_k^\frac 1 k = \tanh(u),\qquad F_k^i = \frac 1 n \delta_k^i\ .\]
At such a $u$, the linearisation becomes
\[0=a^{ij}\tilde \n_{ij} v +b^i \tilde\n_i v +c v\ ,\]
where $a$ is elliptic, $b$ is bounded, and
\begin{flalign*}
c&=\frac{\pard{\tau}{r}}{\tau}H_k^\frac 1 k-3\tanh(u)H_k^\frac 1 k+\tau \frac{\cosh^2(u)+\sinh^2(u)}{\cosh^3(u)}\text{tr}F^{ij}-\psi_r - \psi_\tau \pard \tau z\\
&=\cosh^{-2}(u)-\cosh^p(u)\tanh(u)-u\cosh^{p-2}(u) - p\cosh(u)^{p-1}u\tanh(u) \sinh(u) . 
\end{flalign*}
We recall that $u=\lambda$ was chosen so that $\lambda\cosh^p\lambda=1$, which implies that $u\cosh^{p-2}u=\cosh^{-2}(u)$, and so we see that $c< 0$. The strong maximum principle now implies that the only solution $v\in C^{4,\alpha}(S^n)$ of $\Phi_u(\cdot,0)(v)=0$ is $v=0$. This implies that $\text{ker}(\Phi_u)=\{0\}$ and so, the standard theory of second order elliptic equations imply that $\Phi_u$ is invertable, as required.
\end{proof}

\begin{proof}[Proof of Step 2]
By assumption, we have that $\psi(\xi, z, \cosh(z))<\tanh(z)$ for $z<R_1$ and $\psi(\xi, z, \cosh(z))>\tanh(z)$ for $z>R_2$. Similarly we see directly that there exists $R^\Psi_1, R^\Psi_2>0$ such that 
$\Psi(\xi, z, \cosh(z))<\tanh(z)$ for $z<R^\Psi_1$ and $\psi(\xi, z, \cosh(z))>\tanh(z)$ for $z>R^\Psi_2$. Setting $\overline{R}_1=\min\{R_1,R^\Psi_1\}$ and $\overline{R}_2=\max\{R_2,R^\Psi_2\}$, then for all $t\in[0,1]$, $\psi_t(\xi, z, \cosh(z))<\tanh(z)$ for $z<\overline{R}_1$ and $\psi_t(\xi, z, \cosh(z))>\tanh(z)$ for $z>\overline{R}_2$. Lemma \ref{BarrierConditions}, therefore yields uniform $C^0$ estimates 
\[0<\overline{R}_2\leq u_t\leq \overline{R}_2<\infty . \]
Proposition \ref{TiltEstimate}, in addition to giving a $C^1$ estimate,  implies uniform spacelikeness, and so we may apply Theorem \ref{CurvatureEstimate} to yield $|\lambda_i|<C_A$, which implies uniform $C^2$ estimates in this situation. Uniform parabolicity of the equation now follows, and so due to the classical regularity theory for uniformly elliptic equations and the Evans-Krylov theorem \cite{Evans, Krylov} we obtain
\begin{equation}
  \label{eq:102}
  \|u_t\|_{C^{4,\alpha}(\mathbb{S}^n)} < C,
\end{equation}
for any admissible solution $u_t\in C^{4,\alpha}_{ad}(\mathbb{S}^n)$, where the constant $C$ depends on
$k, n, R_1, R_2$ and $\|\psi\|_{C^{2,\alpha}(\mathbb{S}^n)}$.
\noindent
Due to compactness and the above estimates, there exists a constant $\delta>0$ such that $\delta<\psi_t(\xi, u(\xi), \tau(u))$ for all $\xi\in S^n$. We define the bounded open set $V:=\{\lambda\in\Gamma_k : H_k^\frac 1 k(\lambda)\geq \delta , |\lambda|<\sqrt{n}C_A \}\subset\Gamma_k$ and we define the bounded set
\begin{align*}
  \label{eq:104}
  \mathcal{B} &=  \Big\{u \in C^{4,\alpha}(\mathbb{S}^n) \Big|\\
  &\qquad\qquad\frac 1 2 \overline{R}_1 < u < 2 \overline{R}_2 ,  \| u \|_{C^{4,\alpha}(\mathbb{S}^n)} < C
\mbox{ and } \lambda(A[u(\xi)])\in V \ \forall \ \xi\in\mathbb{S}^n\Big\}.
\end{align*}
\noindent
Clearly the arguments of the previous paragraphs imply that any admissable solution $u_t\in C^{4,\alpha}_{ad}(S^n)$ is contined in $\mathcal{B}$, and $\partial \mathcal{B}\cap \Phi^{-1}(\cdot,t)=\emptyset$ for all $t\in[0,1]$. 
\end{proof}

\begin{proof}[Proof of Step 3]
This step now follows exactly as in \cite{YYLi:elli}. Due to Step 2 and \cite[Definition 2.2, Proposition 2.2]{YYLi:1989}, the degree is defined and constant for $t\in[0,1]$. By \cite[Proposition 2.3, Proposition 2.2]{YYLi:1989}, $\text{deg}(\Phi(\cdot,0), \mathcal{B}, 0)=\text{deg}(\Phi_u(\cdot,0), B_1, 0)$ , where $B_1$ is the open unit ball in $C^{4, \alpha}(S^n)$. However, Step 1 and \cite[Proposition 2.4]{YYLi:1989} imply that $\text{deg}(\Phi(\cdot,0), \mathcal{B}, 0)=\pm 1=\text{deg}(\Phi(\cdot,1), \mathcal{B}, 1)$, and we conclude that a $k$-admissable solution $u_1$ to $\Phi(u_1,1)=0$ exists. Standard elliptic estimates imply that $u_1$ is smooth and therefore the proof of Theorem \ref{FullTheorem} is complete.
\end{proof}


\end{document}